\documentclass[11pt, a4paper]{article}
\usepackage{mathrsfs}
\usepackage{graphicx}
\usepackage{lscape}
\usepackage{enumerate}

\usepackage{amsmath}
\usepackage{amssymb}
\usepackage{amsbsy}

\usepackage{amsfonts}

\newtheorem{lem}{Lemma}[section]%
\newtheorem{theorem}[lem]{Theorem}%
\newtheorem{defi}[lem]{Definition}%
\newtheorem{prop}[lem]{Proposition}%

\def\a{\alpha} \def\b{\beta}   
 \def\s{\sigma}   
 \def\ld{\lambda}

\def\G{\Gamma}

\def\di{\bigm|} \def\lg{\langle} \def\rg{\rangle}
\def\nd{\mathrel{\bigm|\kern-.7em/}}

\def\f{\noindent}

\def\Aut{\hbox{\rm Aut}}

\def\Cay{\hbox{\rm Cay}}

\def\Cos{\hbox{\rm Cos}}
\def\mod{\hbox{\rm mod }}

\def\demo{\f {\bf Proof.}\hskip10pt}

\newcommand{\qed}{\mbox{\raisebox{0.7ex}{\fbox{}}} \vspace{4truemm}}
\def\mz{{\mathbb Z}}
\def\C{{\bf C}}

\def\B{\mathcal{B}}

\def\MG{\mathcal{G}}
\def\K{{K}}
\def\1{{\bf 1}}

\usepackage{color}

\definecolor{Thistle}{rgb}{0.847,0.749,0.847}
\definecolor{Khaki}{rgb}{0.941,0.902,0.549}
\definecolor{Orchid}{rgb}{0.855,0.439,0.839}
\definecolor{MediumOrchid}{rgb}{0.729,0.333,0.827}

\definecolor{brown}{rgb}{0.8,0.5,0}
\definecolor{LightBrown}{rgb}{0.8,0.2,0.4}

\definecolor{DarkGray}{rgb}{0.78,0.78,0.78}
\definecolor{DarkMidGray}{rgb}{0.81,0.81,0.81}
\definecolor{MidGray}{rgb}{0.85,0.85,0.85}
\definecolor{LightGray}{rgb}{0.88,0.88,0.88}
\definecolor{VeryLightGray}{rgb}{0.96,0.96,0.96}

\definecolor{GrayA}{rgb}{0.7,0.7,0.7}
\definecolor{GrayB}{rgb}{0.78,0.78,0.78}
\definecolor{GrayC}{rgb}{0.80,0.80,0.80}
\definecolor{GrayD}{rgb}{0.82,0.82,0.82}
\definecolor{GrayE}{rgb}{0.84,0.84,0.84}
\definecolor{GrayF}{rgb}{0.86,0.86,0.86}
\definecolor{GrayG}{rgb}{0.88,0.88,0.88}
\definecolor{GrayH}{rgb}{0.90,0.90,0.90}
\definecolor{GrayI}{rgb}{0.92,0.92,0.92}
\definecolor{GrayJ}{rgb}{0.94,0.94,0.94}

\definecolor{VeryLightBlue}{rgb}{0.9,0.9,1}
\definecolor{LightBlue}{rgb}{0.8,0.8,1}
\definecolor{MidBlue}{rgb}{0.5,0.5,1}
\definecolor{DarkBlue}{rgb}{0,0,0.6}

\definecolor{Gold}{rgb}{1,0.843,0}

\definecolor{LightGreen}{rgb}{0.88,1,0.88}
\definecolor{MidGreen}{rgb}{0.6,1,0.6}
\definecolor{DarkGreen}{rgb}{0,0.6,0}

\definecolor{VeryLightYellow}{rgb}{1,1,0.9}
\definecolor{LightYellow}{rgb}{1,1,0.6}
\definecolor{MidYellow}{rgb}{1,1,0.5}
\definecolor{DarkYellow}{rgb}{1,1,0.2}

\definecolor{VeryLightRed}{rgb}{1,0.9,0.9}
\definecolor{LightRed}{rgb}{1,0.8,0.8}
\definecolor{MidRed}{rgb}{1,0.55,0.55}

\long\def\delete#1{}

\begin{document}

\title{Weak metacirculants of odd prime power order}
\vspace{4 true mm}


\author {{\sc Jin-Xin Zhou}\\
{\small\em Department of Mathematics}\\
{\small\em Beijing Jiaotong University}\\
{\small\em Beijing 100044, P.R. China}\\
{ \texttt{jxzhou@bjtu.edu.cn}}\\[1ex]
{\sc Sanming Zhou}\\
{\small\em School of Mathematics and Statistics}\\
{\small\em The University of Melbourne}\\
{\small\em Parkville, VIC 3010, Australia}\\
{\texttt{sanming@unimelb.edu.au}} }

\date{}
\openup 0.5\jot

\maketitle

\begin{abstract}
Metacirculants are a basic and well-studied family of vertex-transitive graphs, and weak metacirculants are generalizations of them. A graph is called a weak metacirculant if it has a vertex-transitive metacyclic automorphism group. This paper is devoted to the study of weak metacirculants with odd prime power order. We first prove that a weak metacirculant of odd prime power order is a metacirculant if and only if it has a vertex-transitive split metacyclic automorphism group. We then prove that for any odd prime $p$ and integer $\ell\geq 4$, there exist weak metacirculants of order $p^\ell$ which are Cayley graphs but not Cayley graphs of any metacyclic group; this answers a question in Li et al. \cite{LWS} (2013). We construct such graphs explicitly by introducing a construction which is a generalization of generalized Petersen graphs. Finally, we determine all smallest possible metacirculants of odd prime power order which are Cayley graphs but not Cayley graphs of any metacyclic group.

\medskip

\end{abstract}

\thispagestyle{empty}

\section{Introduction}

Let $m\geq 1$ and $n\geq 2$ be integers. A graph $\G$ of order $mn$ is called \cite{MS} an {\em $(m, n)$-metacirculant graph} (in short {\em $(m, n)$-metacirculant}) if it has an automorphism $\s$ of order $n$ such that $\lg\s\rg$ is semiregular on the vertex set of $\G$, and an automorphism $\tau$ normalizing $\lg\s\rg$ and cyclically permuting the $m$ orbits of $\lg\s\rg$ such that $\tau$ has a cycle of size $m$ in its cycle decomposition. A graph is called a {\em metacirculant} if it is an $(m, n)$-metacirculant for some $m$ and $n$. It follows from this definition that a metacirculant $\G$ has an autormorphism group $\lg \s, \tau\rg$ which is metacyclic and transitive on the vertex set of $\G$. In general, a group $G$ is called {\em metacyclic} if it contains a cyclic normal subgroup $N$ such that $G/N$ is cyclic. In other words, a metacyclic group $G$ is an extension of a cyclic group  $N\cong C_n$ by a cyclic group $G/N\cong C_m$, written $G\cong C_n.C_m$. If this extension is split, namely $G \cong C_n: C_m$, then $G$ is called a {\em split metacyclic group}.

Introduced by Alspach and Parsons \cite{AP}, metacirculants form a basic class of vertex-transitive graphs. As a generalization of metacirculants, Maru\v si\v c and \v Sparl~\cite{MS} introduced the following concept: A graph is called a {\em weak metacirculant} if it has a vertex-transitive metacyclic automorphism group. In \cite{LWS}, Li et al. divided the class of weak metacirculants into the following two subclasses:  A weak metacirculant is called a {\em split weak metacirculant} or {\em non-split weak metacirculant} according to whether or not its full automorphism group contains a vertex-transitive split metacyclic subgroup. In \cite{LWS}, Li et al. studied the relationship between metacirculants and weak metacirculants. Among other results they proved that every metacirculant is a split weak metacirculant (see \cite[Lemma~2.2]{LWS}), but it was unknown whether the converse of this statement is true. So the following question arises naturally.

\medskip
\f{\bf Question~A}\, Is it true that any split weak metacirculant is a metacirculant?
\medskip

In this paper we first give a positive answer to this question for split weak metacirculants of odd prime power order, as stated in the following result.

\begin{theorem}\label{p-metacirculant}
A connected weak metacirculant with order an odd prime power is a metacircualnt if and only if it is a split weak metacirculant.
\end{theorem}

Question A is open for split weak metacirculants of order not an odd prime power; in fact, there is no result concerning Question A in the literature in this case as far as we know.

Obviously, any Cayley graph of a metacyclic group is a weak metacirculant; such a graph is called a {\em weak metacirculant Cayley graph} (see \cite[p.41]{LWS}). Weak metacirculant Cayley graphs form a large class of weak metacirculants. However, not every weak metacirculant is a Cayley graph. For example, the Petersen graph is a $(2, 5)$-metacirculant but not a Cayley graph. The following question was posed by Pan \cite[p.15]{Pan-thesis} and Li et al.~\cite[p.41]{LWS} independently.

\medskip
\f{\bf Question~B}\, Is it true that a weak metacirculant which is a Cayley graph of some (not necessarily metacyclic) group must be a weak metacirculant Cayley graph?
\medskip

Our second main result gives a negative answer to this question.

\begin{theorem}\label{exist}
Let $p$ be an odd prime. Then for any integer $\ell\geq 4$ there exists a weak metacirculant of order $p^\ell$ which is a Cayley graph but not a weak metacirculant Cayley graph.

Moreover, the smallest possible order and valency of a weak metacirculant with order a power of $p$ which is a Cayley graph but not a weak metacirculant Cayley graph are $p^4$ and $2p+2$, respectively.
\end{theorem}

The third main result in this paper is the following classification of connected metacirculants of order $p^4$ and valency $2p+2$, where $p$ is an odd prime. The graph ${\rm MP}_{p^3, p^2, p^2, \ld}$ involved in the classification will be defined in Definition~\ref{multi}; it belongs to a large family of graphs that contains all generalized Petersen graphs as a proper subfamily.

\begin{theorem}
\label{2p+2}
Let $p$ be an odd prime. Let $\G$ be a connected metacirculant of order $p^4$ and valency $2p+2$. Then one of the following holds:
\begin{enumerate}[\rm (a)]
  \item $\G$ is a metacirculant Cayley graph;
  \item $\G$ is not a Cayley graph;
  \item $\G$ is isomorphic to ${\rm MP}_{p^3, p^2, p^2, \ld}$ for some element $\ld$ of $\mz_{p^3}^*$ with order $p^2$.
\end{enumerate}
\end{theorem}

This result seems to suggest that most weak metacirculants which are Cayley graphs are weak metacirculant Cayley graphs. Nevertheless, more research is needed to find out whether this is indeed the case.

The rest of this paper is organized as follows. In the next section we will collect some basic definitions on permutation groups, Cayley graphs and vertex-transitive graphs. In section \ref{sec:p-meta}, we will give the proof of Theorem \ref{p-metacirculant} after presenting a few results on $p$-groups. In section \ref{sec:weak meta Cay}, we will prove that any weak metacirculant of order an odd prime power $p^{n}$ must be a weak metacirculant Cayley graph if its valency is less than $2p+2$ or its order is at most $p^3$. This result will be used in the proof of Theorem \ref{exist}, which will be given in section \ref{sec:exist}. Another preparation for the proof of Theorem \ref{exist} is the construction of multilayer generalized Petersen graphs, which will be introduced in section \ref{sec:multi Petersen}. The proof of Theorem \ref{2p+2} will be given in section \ref{sec:2p+2}.

\section{Preliminaries}

\subsection{Definitions and notation}

Given a group $G$, denote by $1_G$, $\Aut(G)$, $Z(G)$, $\Phi(G)$ and $G'$ the identity element, full automorphism group, center, Frattini subgroup and derived subgroup of $G$, respectively. Denote by $o(x)$ the order of an element $x$ of $G$. For a subgroup $H$ of $G$, denote by $C_G(H), N_G(H)$ the centralizer and normalizer of $H$ in $G$, respectively. Of course $C_G(H)$ is normal in $N_G(H)$, and the well-known N/C theorem asserts that the quotient group $N_G(H)/C_G(H)$ is isomorphic to a subgroup of $\Aut(H)$. Given a $p$-group $G$ of exponent $p^e$, where $p$ is a prime and $e \ge 1$ an integer, for each integer $s$ between $0$ and $e$, set
$$
\Omega_s(G)=\lg g\in G\ |\ g^{p^s}=1_G\rg.
$$

A {\em block of imprimitivity} of a permutation group $G$ on a set $\Omega$
is a subset $\Delta$ of $\Omega$ with $1 < |\Delta| < |\Omega|$ such that for any $g \in G$,
either $\Delta^g=\Delta$ or $\Delta^g \cap \Delta=\emptyset$. In this case the {\em blocks} $\Delta^g,\ g \in G$ form a {\em $G$-invariant partition} of $\Omega$.

We reserve $C_n$ for the cyclic group of order $n$, $\mz_n$ for the ring of integers modulo $n$, and $\mathbb{Z}_{n}^{*}$ for the multiplicative group of units of $\mathbb{Z}_{n}$ consisting of integers coprime to $n$.

All graphs in this paper are finite, simple and undirected. For a graph $\G$, we denote its vertex set and edge set by $V(\G)$ and $E(\G)$, respectively. Given $u, v\in V(\G)$, denote by $u \sim v$ the relation that $u$ is adjacent to $v$ in $\G$, by $\{u, v\}$ the edge between $u$ and $v$, and by $(u, v)$ the arc from $u$ to $v$. Denote by $\G(v)$ the neighbourhood of $v$, and by $\G[B]$ the subgraph of $\G$ induced by a subset $B$ of $V(\G)$. An {\em $s$-cycle} in $\G$, denoted by $\C_s$, is an $(s+1)$-tuple of  pairwise distinct vertices $(v_0, v_1, \ldots, v_s)$ such that $\{v_{i-1 }, v_i\}\in E(\G)$ for $1\leq i\leq s$ and $\{v_s, v_0\}\in E(\G)$. Denote by $\K_n$ the complete graph of order $n$, and $\K_{n, n}$ the complete bipartite graph with biparts of cardinality $n$. The {\em lexicographic product} of a graph $\G_1$ by a graph $\G_2$, denoted by $\G_1\circ\G_2$, is the graph with vertex set $V(\G_1)\times V(\G_2)$ such that $\{(x_1, x_2), (y_1, y_2)\}\in E(\G_1\circ\G_2)$ if and only if either $\{x_1, y_1\}\in E(\G_1)$, or $x_1=y_1$ and $\{x_2, y_2\}\in E(\G_2)$.

The full automorphism group of a graph $\G$ is denoted by $\Aut(\G)$. $\G$ is called {\em $G$-vertex-transitive} (respectively, {\em $G$-edge-transitive}) if $G \le \Aut(\G)$ and $G$ is transitive on $V(\G)$ (respectively, $E(\G)$); in this case $G$ is said to be a vertex-transitive (respectively, edge-transitive) automorphism group of $\G$. $\G$ is {\em vertex-transitive} (respectively, {\em edge-transitive}) if it is $\Aut(\G)$-vertex-transitive (respectively, $\Aut(\G)$-edge-transitive). {\em $G$-arc-transitive} graphs and {\em arc-transitive} graphs are understood similarly. Given a $G$-vertex-transitive graph $\G$ and a $G$-invariant partition $\B$ of $V(\G)$, the {\em quotient graph} of $\G$ with respect to $\B$, denoted by $\G_\B$, is defined as the graph with vertex set $\B$ such that, for distinct $B, C\in \B$, $B$ is adjacent to $C$ if and only if
there exist $u\in B$ and $v\in C$ which are adjacent in $\G$. In particular, for a normal subgroup $N$ of $G$, the set $\B$ of orbits of $N$ on $V(\G)$ is a $G$-invariant partition of $V(\G)$, and in this case we use $\G_N$ in place of $\G_\B$.

\subsection{Cayley graphs}

Given a finite group $G$ and an inverse-closed subset $S\subseteq G\setminus\{1_G\}$, the {\em Cayley graph} $\Cay(G,S)$ of $G$ with respect to $S$ is the graph with vertex set $G$ and edge set $\{\{g,sg\}\mid g\in G,s\in S\}$. It is well known that the right regular representation $R(G)=\{R(g)\ |\ g\in G\}$ of $G$ is a subgroup of $\Aut(\Cay(G,S))$, where $R(g)$ is the permutation of $G$ defined by $R(g): x\mapsto xg$ for $x\in G$. In \cite{Godsil1981}, Godsil proved that the normalizer of $R(G)$ in $\Aut(\Cay(G,S))$ is $R(G): \Aut(G,S)$, where $\Aut(G,S)$ is the group of automorphisms of $G$ fixing $S$ setwise. In the case when $R(G)$ is normal in $\Aut(\Cay(G,S))$, $\Cay(G,S)$ is called \cite{X1} a {\em normal Cayley graph}. The reader is referred to \cite{FengLuXu} for recent results on normal Cayley graphs.

It is well known that a graph $\G$ is isomorphic to a Cayley graph if and only if it has an automorphism group acting regularly on its vertex set (see \cite[Lemma~16.3]{B}). In general, a permutation group $G$ on a set $\Omega$ is called {\em semiregular} on $\Omega$ if $G_\a=1_G$ for every $\a\in \Omega$, and {\em regular} on $\Omega$ if $G$ is transitive and semiregular on $\Omega$, where $G_\a$ is the \emph{stabilizer} of $\a$ in $G$, defined as the subgroup of $G$ consisting of those elements of $G$ which fix $\a$.

\subsection{Coset graphs}

Let $G$ be a finite group, $H$ a subgroup of $G$, and $D$ the union of
some double-cosets $HgH$ with $g\notin H$ such that
$D=D^{-1}$. The {\em coset graph} $\G=\Cos(G, H, D)$ of $G$ with
respect to $H$ and $D$ is defined as the graph with vertex set
$V(\G)=[G: H]$, the set of right cosets of $H$ in $G$, and edge set
$E(\G)=$ $\{\{Hg, Hdg\}\di g\in G,$ $ d\in D \}$. It is easy to see
that $\G$ is well defined and has valency $|D|/|H|$. Further, $\G$ is
connected if and only if $D$ generates $G$. In the special case when $H=1_G$,
$\G$ is the Cayley graph of $G$ with respect to $D$.
Denote by $R_H$ the right multiplication action of $G$ on
$V(\G)=[G: H]$, defined by $R_H(g): Hx \mapsto Hxg,\ Hx \in [G:H]$. (In particular, $R_{1_G}(G)$ is the right regular representation
$R(G)$ of $G$.) Then $R_H$ is transitive on $V(\G)$, and $R_H$ is faithful on $V(\G)$ if and only if $H$ is \emph{core-free} in $G$, that is, $\cap_{g\in G}H^g=1_G$. It is easy to see that $R_H(G) \leq \Aut(\G)$. Hence $\G$ is
vertex-transitive. In \cite{Sabidussi}, Sabidussi proved that all vertex-transitive graphs can be constructed this way up to isomorphism.

\begin{prop}\label{coset}
The coset graph $\Cos(G, H, D)$ constructed above is $G$-vertex-transitive. Conversely, if $\G$ is a $G$-vertex-transitive graph, then it is isomorphic to a
coset graph $\Cos(G, H, D)$, where $H=G_\a$ for a fixed $\a \in V(\G)$ and $D$ consists of all elements of $G$ which map $\a$
to one of its neighbours.
\end{prop}

The following results are well known; see, for example, \cite[Lemma~2.1]{LLZ}.

\begin{lem}
\label{lem-coset}
Let $G$ be a group and $H$ a core-free subgroup of $G$. Take $g\in G \setminus H$ and let $\G = \Cos(G, H, H\{g,g^{-1}\}H)$. Then the following hold:
\begin{enumerate}[\rm (a)]
\item $\G$ is $G$-edge-transitive;
\item $\G$ is $G$-arc-transitive if and only if $HgH=Hg^{-1}H$;
\item $\G$ is connected if and only if $G=\lg H, g\rg$;
\item the valency of $\G$ is equal to $|H:H\cap H^g|$ if $HgH=Hg^{-1}H$, or $2|H: H\cap H^g|$ if $HgH \ne Hg^{-1}H$.
\end{enumerate}
\end{lem}

\section{Proof of Theorem~\ref{p-metacirculant}}
\label{sec:p-meta}

\subsection{Some results on $p$-groups}

In order to prove Theorem~\ref{p-metacirculant}, we first present a few results on $p$-groups. The following result is due to Xu and Zhang; see \cite[Theorem~2.1]{Xu-Zhang}.

\begin{lem}
\label{p-metacyclic}
Let $p$ be an odd prime and $G$ a metacyclic $p$-group. Then $G$ has representation
$$G=\lg a, b \ |\ a^{p^{r+s+u}}=1_G,\ b^{p^{r+s+t}}=a^{p^{r+s}},\ b^{-1}ab=a^{1+p^r}\rg$$
for some nonnegative integers $r, s, t, u$ such that $r\geq 1, r\geq u$. Moreover, different values of the parameters $r,s,t,u$ satisfying these conditions give rise to non-isomorphic metacyclic $p$-groups, and $G$ is non-split if and only if $stu\neq 0$.
\end{lem}

The following result can be easily proved (see, for example, \cite[Exercise 85]{p-group1}).

\begin{lem}
\label{p2-group}
Let $G$ be a noncyclic metacyclic $p$-group. If $p>2$, then $\Omega_1(G)$ is a normal subgroup of $G$ which is isomorphic to $C_p\times C_p$.
\end{lem}

A $p$-group $G$ is said to be {\em $p^r$-abelian} if $(xy)^{p^r}=x^{p^r}y^{p^r}$ for any $x, y\in G$.

\begin{lem}\label{p-abelian}
Any metacyclic $p$-group $G$ with $p>2$ is $p^\ell$-abelian, where $p^\ell=|G'|$.
\end{lem}

\demo By \cite[Theorem~7.1 (c)]{p-group1}, $G$ is regular, and then by \cite[III, 10.8(g)]{Huppert}, $G$ is $p^\ell$-abelian.
\hfill\qed

\begin{lem}\label{complement}
Let $p$ be an odd prime. Let $G=\lg \s\rg:\lg\tau\rg\cong C_{p^m}: C_{p^n}$ with $m\geq n\geq 1$. For any $1_G\neq g\in G$, if $\lg g\rg\cap\lg\s\rg=1_G$, then there exists $\tau'\in G$ such that $\lg\tau'\rg\cong C_{p^n}$ and $g\in\lg\tau'\rg$.
\end{lem}

\demo We make induction on the order $|G|$ of $G$. Clearly, $|G|=p^{m+n}\geq p^2$. Suppose that $|G|=p^2$. Then $m=n=1$ and $G=\lg \s\rg:\lg\tau\rg\cong C_{p}\times C_{p}$. For any $1_G\neq g\in G$, if $\lg g\rg\cap\lg\s\rg=1_G$, then $\lg g\rg\cong C_p$ and $G=\lg\s\rg\times\lg g\rg$, as required.

In what follows we assume that $|G|>p^2$. Let $1_G\neq g\in G$ be such that $\lg g\rg\cap\lg\s\rg=1_G$. Since $G/\lg \s\rg\cong C_{p^n}$, $g$ has order at most $p^n$. If the order $o(g)$ of $g$ is equal to $p^n$, then the result is clearly true. Assume $o(g)=p^k<p^n$. Then $n>k\geq 1$.

Suppose first that $m>n$. Then $g, \tau\in\Omega_{m-1}(G)$. Since $G/\lg\s\rg\cong C_{p^n}$, we have $G'\leq\lg\s\rg$. Clearly, $G'\neq\lg\s\rg$, so $G'\leq\lg\s^p\rg\cong\mz_{p^{m-1}}$. By Lemma~\ref{p-abelian}, $G$ is $p^{m-1}$-abelian. This implies that $\Omega_{m-1}(G)$ contains no elements of order greater than $p^{m-1}$, and consequently, $\Omega_{m-1}(G)<G$. Clearly, $\s^p\in \Omega_{m-1}(G)$, so $\Omega_{m-1}(G)=\lg \s^p\rg: \lg \tau\rg\cong C_{p^{m-1}}: C_{p^n}$. Since $m-1\geq n$, by induction there exists $\tau'\in\Omega_{m-1}(G)$ such that $\lg\tau'\rg\cong C_{p^n}$ and $g\in\lg\tau'\rg$, as required.

Now suppose that $m=n$. Since $p>2$, $\Aut(C_{p^m})\cong C_{p^{m-1}(p-1)}$. Then $\tau$ induces an automorphism by conjugation of $\lg\s\rg$ of order at most $p^{m-1}$. This implies that $\tau^{p^{m-1}}$ commutes with $\s$, and so $\tau^{p^{m-1}}$ is in the center of $G$. Since $p>2$, by Lemma~\ref{p2-group} we have $\Omega_1(G)\cong C_p\times C_p$ and so $\Omega_1(G)=\lg\s^{p^{m-1}}\rg\times\lg\tau^{p^{m-1}}\rg\cong C_{p}\times C_p$.

Let $N=\lg\tau^{p^{m-1}}\rg$. Then $G/N\cong C_{p^m}: C_{p^{m-1}}$. If $gN=N$, then $g\in N$ and $g\in\lg\tau\rg$, and so the result holds.

Assume that $gN\neq N$ in the sequel. If $\lg gN\rg\cap\lg\s N\rg=N$, then by induction, $gN\in\lg\tau'N\rg$ for some $\tau'\in G$ such that $\lg\tau'N\rg\cap\lg\s N\rg=N$ and $\lg\tau'N\rg\cong C_{p^{m-1}}$. So $g\in \lg\tau', N\rg$. If $\lg\tau'\rg\cap N=1_G$, then $\Omega_1(G)\leq\lg\tau', N\rg$ and so $\lg\s^{p^{m-1}}N\rg=\Omega_1(G)/N\leq \lg\tau', N\rg/N=\lg\tau' N\rg$, a contradiction. Thus $\lg\tau'\rg\cap N\neq1_G$, and hence $N\leq \lg\tau'\rg$. This implies that $g\in\lg\tau'\rg\cong C_{p^n}$, as required.

Suppose that $\lg gN\rg\cap\lg\s N\rg\neq N$. Then $g^iN=\s^j N\neq N$ for some $i,j$. It then follows that $g^i=\s^j (\tau^{p^{m-1}})^k$ for some $k$. If $p\ |\ k$, then $g^i=\s^j\neq 1_G$, a contradiction. Thus, $(p, k)=1$. Let $j=p^{\ell}j'$ be such that $(j', p)=1$. Then $\ell\leq m-1$. If $\ell<m-1$, then $\s^j$ has order at least $p^2$ and so $g^{pi}=\s^{pj}=\s^{p^{\ell+1}j'}\neq1_G$, which contradicts the fact that $\lg g\rg\cap\lg\s\rg=1_G$. Hence, $\ell=m-1$, and so $g^i=(\s^{j'}\tau^{k})^{p^{m-1}}$ as $G$ is $p^{m-1}$-abelian.

Since $(k, p)=1$, we have $G=\lg \s, \s^{j'}\tau^k\rg$, and since $1\neq g^i\in\lg \s^{j'}\tau^k\rg$ and $\lg\s\rg\cap\lg g\rg=1_G$, we have $\lg\s\rg\cap\lg \s^{j'}\tau^k\rg=1_G$. This implies that $G=\lg \s\rg: \lg\s^{j'}\tau^k\rg\cong C_{p^m}: C_{p^m}$ and so $\lg\s^{j'}\tau^{k}\rg\cong C_{p^m}$. If $g$ has order $p$, then $\lg g\rg=\lg g^i\rg\leq \lg \s^{j'}\tau^{k}\rg$, as required. If $g$ has order $p^t$ with $t>1$, then let $M=\lg g^{p^{t-1}}\rg=\lg(\s^{j'}\tau^{k})^{p^{m-1}}\rg$. Clearly, $G/M\cong C_{p^{m}}:C_{p^{m-1}}$ and $\lg gM\rg\cap\lg\s M\rg=M$. By induction, $gM\in\lg \tau'M\rg$  for some $\tau'\in G$ such that $\lg\tau' M\rg\cap\lg\s M\rg=M$ and $\lg\tau'M\rg\cong C_{p^{m-1}}$. If $\lg\tau'\rg\cap M=1_G$, then $\Omega_1(G)\leq\lg\tau', M\rg$ and so $\lg\s^{p^{m-1}}M\rg=\Omega_1(G)/M\leq \lg\tau', M\rg/M=\lg\tau' M\rg$, a contradiction. Thus, $\lg\tau'\rg\cap M\neq 1_G$, and hence $M\leq \lg\tau'\rg$. This implies that $g\in\lg\tau'\rg\cong C_{p^m}$, as required.\hfill\qed

\subsection{Proof of Theorem~\ref{p-metacirculant}}

We prove the following result first.

\begin{lem}\label{cyclic-center}
Let $\G$ be a connected weak metacirculant with order a power of an odd prime $p$. Then $\Aut(\G)$ contains a metacyclic $p$-subgroup $G$ which is transitive on $V(\G)$. Moreover, if $Z(G)$ is not cyclic, then $G$ is regular on $V(\G)$ and so $\G$ is a weak metacirculant Cayley graph.
\end{lem}

\demo Since $\Gamma$ is a weak metacirculant, $\Aut(\G)$ has a metacyclic subgroup $X$ which is transitive on $V(\G)$. Let $G$ be a Sylow $p$-subgroup of $X$. Then $G$ is metacyclic, and by \cite[Theorem~3.4]{WI}, $G$ is also transitive on $V(\G)$, proving the first statement in the lemma.

Since $p>2$, we have $\Omega_{1}(G)\cong C_p\times C_p$ by Lemma~\ref{p2-group}. If $Z(G)$ is not cyclic, then $\Omega_1(G)\leq Z(G)$. For any $v\in V(\G)$, if $G_v\neq 1_G$, then $G_v\cap \Omega_1(G)\neq 1_G$. However, $G_v\cap \Omega_1(G)\unlhd G$. So $G_v\cap \Omega_1(G)$ fixes every vertex of $\G$, a contradiction. Thus, $G_v=1_G$, and so $G$ is regular on $V(\G)$. It follows that $\G$ is a Cayley graph of $G$. \hfill\qed

Now we are ready to prove Theorem~\ref{p-metacirculant}.\medskip


\f{\bf Proof of Theorem~\ref{p-metacirculant}.}\; By \cite[Lemma~2.2]{LWS}, each metacirculant has a vertex-transitive split metacyclic automorphism group. The necessity follows. It remains to prove the sufficiency.

Suppose that $G$ is a split metacyclic vertex-transitive $p$-subgroup of $\Aut(\G)$. If $G$ is regular on $V(\G)$, then $\G$ is a Cayley graph of $G$. Since $G$ is a split metacyclic group, $\G$ is a metacirculant graph, as required. In what follows we assume that $G$ is not regular on $V(\G)$. Then $Z(G)$ must be cyclic by Lemma~\ref{cyclic-center}.

\medskip
\noindent {\bf Claim.}\; $G$ can be written as $\lg x\rg:\lg y\rg\cong C_{p^m}: C_{p^n}$ for some integers $m\geq n$.\medskip

In fact, by Lemma~\ref{p-metacyclic}, we have
$$
G=\lg a, b \ |\ a^{p^{r+s+u}}=1_G,\ b^{p^{r+s+t}}=a^{p^{r+s}},\ b^{-1}ab=a^{1+p^r}\rg
$$
for some nonnegative integers $r, s, t, u$ with $r \ge 1, r \ge u$. A straightforward computation leads to the following observations:
\begin{enumerate}[(i)]
   \item $|G|=p^{2(r+s)+u+t}$, {\rm exp}$(G)=p^{r+s+t+u}=o(b)$;
  \item $G'=\lg a^{p^r}\rg\cong C_{p^{s+u}}$ and $G$ is $p^{s+u}$-abelian;
  \item $Z(G)=\lg a^{p^{s+u}}, b^{p^{s+u}}\rg$.
\end{enumerate}

Since $o(a^{p^{s+u}})=p^r\leq o(b^{p^{s+u}})=p^{r+t}$ and $Z(G)$ is cyclic, we then have $Z(G)=\lg b^{p^{s+u}}\rg$. Consequently, $a^{p^{s+u}}\in \lg a\rg\cap\lg b\rg=\lg b^{p^{r+s+t}}\rg=\lg a^{p^{r+s}}\rg\cong C_{p^u}$ and so $r\leq u$. This together with $u\leq r$ implies $r=u$.

Since $G$ is split, by Lemma~\ref{p-metacyclic}, we have $stu=0$. If $u=0$, then $u=r=0$, contradicting the assumption that $r\geq1$. So $u>0$. We then have
$$
G=\left\{
      \begin{array}{ll}
        \lg b\rg: \lg ab^{-p^t}\rg\cong C_{p^{2r+t}}: C_{p^r}, & \hbox{{\rm if}}\ s=0; \\
        \lg a\rg: \lg ba^{-1}\rg\cong C_{p^{2r+s}}: C_{p^{r+s}}, & \hbox{{\rm if}}\ t=0
      \end{array}
    \right.
$$
as stated in the Claim.

By the Claim above, $G=\lg x\rg:\lg y\rg\cong C_{p^m}: C_{p^n}$ with $m\geq n$. Since $G$ is transitive on $V(\G)$, $\lg x\rg$ acts semiregularly on $V(\G)$. Assume that $\lg x\rg$ has $p^\ell$ orbits for some $\ell\leq n$. For any $v\in V(\G)$, let $G_v=\lg z\rg$. Then $y^{p^\ell}\in \lg x\rg: G_v$ and $|G_v|=p^{n-\ell}$. Moreover, $G_v\cap \lg x\rg=1_G$. It follows that $G_v\cong G_v\lg x\rg/\lg x\rg\leq G/\lg x\rg\cong C_{p^n}$. By Lemma~\ref{complement}, there exists $y'\in G$ such that $\lg y'\rg\cong C_{p^n}$ and $z\in\lg y'\rg$. So $G_v$ is a subgroup of $\lg y'\rg$ of order $p^{n-\ell}$, and hence $\lg (y')^{p^\ell}\rg=G_v$. Since $\lg x\rg\cap G_v=1_G$, we have $\lg x\rg\cap \lg y'\rg=1_G$, and since $\lg y'\rg\cong C_{p^n}$, we have $G=\lg x\rg: \lg y'\rg$. Then $y'$ cyclically permutes the $p^\ell$  orbits of $\lg x\rg$, and $(y')^{p^\ell}\in G_v$, implying that $\G$ is a metacirculant.
\hfill\qed

\section{Smallest possible order and valency}
\label{sec:weak meta Cay}

The main result in this section is the following lemma, which asserts that for an odd prime $p$, if a weak metacirculant of order $p^n$ that is a Cayley graph but not a weak metacirculant Cayley graph exists, then it has order at least $p^4$ and valency at least $2p+2$. In the next two sections we will see that both $p^4$ and $2p+2$ are attainable, as needed to establish the second statement in Theorem \ref{exist}.

\begin{lem}
\label{cay-meta}
Let $p$ be an odd prime. Let $\G$ be a weak metacirculant of order $p^n$ for some integer $n \ge 1$. If $\G$ has valency less than $2p+2$ or $n$ is at most $3$, then $\G$ must be a weak metacirculant Cayley graph.
\end{lem}

\demo By Lemma~\ref{cyclic-center}, $\Aut(\G)$ has a metacyclic $p$-subgroup $G$ which is transitive on $V(\G)$. If $G$ is regular on $V(\G)$, then obviously $\G$ is a weak metacirculant Cayley graph. In what follows we assume that $G$ is not regular on $V(\G)$. Then $G$ is non-abelian, and by Lemma~\ref{cyclic-center}, $Z(G)$ is cyclic. Since $\G$ has odd order $p^n$, its valency must be even. We are going to show that $\G$ is a circulant if it has valency less than $2p+2$ or $1 \le n \le 3$.

If $\G$ has valency less than $2p$, then by \cite[Lemma 2.4]{Feng}, $G$ is regular, which contradicts our assumption. Suppose that $\G$ has valency $2p$. Since $G$ is a metacyclic $p$-group with $p>2$, we have $\Omega_1(G)\cong C_p\times C_p$ by Lemma~\ref{p2-group}. Since $G$ is not regular on $V(\G)$, we have $G_v>1$ for $v\in V(\G)$, and so $G_v\cap\Omega_1(G)>1_G$. Consider the quotient graph $\G_{\Omega_1(G)}$ of $\G$ relative to $\Omega_1(G)$. Each orbit of $\Omega_1(G)$ has length $p$, and the subgraph of $\G$ induced by any two adjacent orbits of $\Omega_1(G)$ is isomorphic to $\K_{p, p}$. So $\G_{\Omega_1(G)}\cong\C_{p^\ell}$ for some integer $\ell\geq 1$. Therefore, $\G\cong \C_{p^\ell}\circ p\K_1$, which is a circulant.

Suppose that $1 \le n \le 3$. Since $G$ is metacyclic, we may assume that $G=\lg\s, \tau\rg$ with $\lg\s\rg\unlhd G$. 
Recall that $G$ is transitive but not regular on $V(\G)$ and $G$ is non-abelian.
Since $\lg\s\rg\unlhd G$, $\lg\s\rg$ is semiregular on $V(\G)$.
In the following we will prove that $\lg\s\rg$ is transitive on $V(\G)$. Once this is achieved, it then follows that $\lg\s\rg$ is regular on $V(\G)$ and so $\G$ is a Cayley graph of $\lg\s\rg$, as required.

Suppose to the contrary that $\lg\s\rg$ is intransitive on $V(\G)$. Since $n\leq 3$, we have $\lg\s\rg\cong C_p$ or $C_{p^2}$. If $\lg\s\rg\cong C_p$, then it is in the center of $G$, and so $G$ is abelian, a contradiction. Thus $\lg\s\rg\cong C_{p^2}$. So $n=3$ and $\tau$ induces an automorphism of $\lg\s\rg$ of order $p$. It follows that $\s^\tau=\s^{kp+1}$ for some integer $1\leq k\leq p-1$. This implies that $G'=\lg\s^p\rg\cong\mz_p$. Since $G$ is a $2$-generator group, by elementary $p$-group theory (see, for example, \cite[Lemma~65.2]{p-group2}), $G$ is an inner abelian $p$-group and therefore $\Phi(G)=Z(G)$. (A group is inner abelian if it is non-abelian but all its proper subgroups are abelian.) Moreover, by \cite{Redei} or \cite[Lemma~65.1]{p-group2}, we may assume that
$$
G=\lg\s,\tau\ |\ \s^{p^2}=\tau^{p^s}=1_G,\ \s^\tau=\s^{p+1}\rg (s\geq1).
$$
Take $v\in V(\G)$. Since $\Phi(G)=Z(G)$, we have $G_v\cap\Phi(G)=1_G$. Since $G$ is not regular on $V(\G)$, we have $G_v\neq1_G$. Take $1_G\neq x\in G_v$. Then $x\notin \Phi(G)$. Since $x^p\in\Phi(G)$, we have $x^p=1_G$. Since $x\notin \Phi(G)=Z(G)$, $x$ commutes with at most one of $\s$ and $\tau$. If $x\s\neq \s x$, then $G=\lg\s, x\rg$ because $G$ is inner-abelian and so $|G|=p^3$. However, this is impossible as $G$ is not regular on $V(\G)$. If $x\s=\s x$, then $x\tau\neq \tau x$ and so $G=\lg\tau, x\rg$. We may write $x=\s^i\tau^j$ for some integers $i, j$. Then $i\in\mz_{p^2}^*$ as $G=\lg\tau, x\rg=\lg\tau, \s^i\tau^j\rg$. Since $p>2$, by Lemma~\ref{p-abelian}, $G$ is $p$-abelian, and hence $(\s^i\tau^j)^p=\s^{pi}\tau^{pj}$. It then follows that $1_G=x^p=\s^{pi}\tau^{pj}$. Thus, $\s^{pi}=1_G$, but this is a contradiction as $\s$ has order $p^2$.
\hfill\qed

\section{Multilayer generalized Petersen graphs}
\label{sec:multi Petersen}

In this section we introduce a construction that can be viewed as a generalization of generalized Petersen graphs. In the next section, we will see that in a special case this construction gives rise to an infinite family of weak metacirculants of odd prime power order which are Cayley graphs but not weak metacirculant Cayley graphs, as needed to establish Theorem \ref{exist}. Introduced in \cite{W}, generalized Petersen graphs are well studied, and they have been generalized in several ways in recent years (see, for example, \cite{CPZ,SPP}). Our generalization is different from the existing ones.


Let $n\geq 3$ and $1\leq t< n/2$. The {\em generalized Petersen graph ${\rm P}(n,t)$} is the graph with vertex set $\{x_i,y_i\ |\ i\in \mz_n\}$ and edge set the union of the {\em out edges} $\{\{x_i,x_{i+1}\}\ |\ i\in\mz_n\}$, the {\em inner edges} $\{\{y_i,y_{i+t}\}\ |\ i\in\mz_n\}$ and the {\em spokes} $\{\{x_i,y_{i}\}\ |\ i\in \mz_n\}$. It is evident that ${\rm P}(n, t)$ has an automorphism $\a=(x_0\ x_1\ x_2\ \ldots\ x_{n-1})(y_0\ y_1\ y_2\ \ldots\ y_{n-1})$ and that $H=\lg\a\rg$ is semiregular on the vertex set of ${\rm P}(n, t)$ with two orbits, namely $X=\{x_i\ |\ i\in\mz_n\}$ and $Y=\{y_i\ |\ i\in\mz_n\}$. The subgraph of ${\rm P}(n, t)$ induced by $X$ is an $n$-cycle while the subgraph of ${\rm P}(n, t)$ induced by $Y$ is the union of some vertex-disjoint cycles.

We now generalize generalized Petersen graphs in such a way that the cyclic semiregular subgroup $H$ has $m\geq 2$ orbits on the vertex set and that the subgraph induced on each orbit of $H$ is a lexicographic product of the union of some cycles of equal length and an empty graph.

\begin{defi}\label{multi}
{\rm
Let $m, n, s$ and $t$ be positive integers such that $m\geq 3, n\geq 2, s\ |\ m$, $1\leq t< \frac{m}{2}$ and $(t, m)=1$. Let $H=\lg h\rg\cong C_m$. For each $i\in\mz_n$, let
$$
\begin{array}{lll}
V_i&=&\{(h^j, i)\ |\ j\in\mz_{m}\},\\
E_i&=&\{\{(h^j, i), (h^{j+ks+t^i}, i)\}, \{(h^j, i), (h^{j-ks-t^i}, i)\}\ |\ k\in\mz_{\frac{m}{s}}, j\in\mz_{m}\},\\
E_{i, i+1}&=&\{\{(h^j, i), (h^j, {i+1})\}\ |\ j\in\mz_{m}\},
\end{array}
$$
where the subscripts are modulo $n$. Define the graph $\G = {\rm MP}_{m, n, s, t}$ by
$$
V(\G) = \cup_{i=0}^{n-1}V_i,\; E(\G) = \cup_{i=0}^{n-1}(E_i \cup E_{i, i+1}), 
$$
with subscripts modulo $n$, and call it the {\em multilayer generalized Petersen graph} with parameters $(m, n, s, t)$.}
\end{defi}

It can be verified that ${\rm MP}_{m, 2, m, t}$ is exactly the generalized Petersen graph ${\rm P}(m, t)$.

For each $g \in H$, define the permutation $R(g)$ on the vertices of $\G = {\rm MP}_{m, n, s, t}$ by
$$
(h^j, i)^{R(g)}=(h^jg, i),\; \mbox{for}\; i \in \mz_n, j\in \mz_m.
$$
Let $R(H)=\{R(g)\ |\ g\in H\}$. A simple computation shows that $R(H)$ is a semiregular subgroup of $\Aut(\G)$ isomorphic to $H$ whose orbits on $V(\G)$ are $V_i$, $i\in\mz_n$. Denote by $\G[V_i]$ the subgraph of $\G$ induced by $V_i$ for $i\in\mz_n$.
Then for each $i \in \mz_n$ the edges between $V_i$ and $V_{i+1}$ form a perfect matching, the subgraph $\G[V_0]$ is isomorphic to the lexicographic product $\C_{s}\circ \frac{m}{s}\K_1$, and for $i\in\mz_n \setminus \{0\}$ the subgraph $\G[V_i]$ is isomorphic to the lexicographic product of the union of some $\ell$-cycles and $\frac{m}{s}\K_1$, where $\ell\ |\ s$. Hence if $n>2$, then $\G$ has valency $2\left(\frac{m}{s}+1\right)$, while if $n=2$, then it has valency $\frac{2m}{s}+1$. 

\begin{lem}\label{block}
Let $\G={\rm MP}_{m, n, s, t}$. If $m/s>1$, then each $V_i$ is a block of imprimitivity of $\Aut(\G)$ on $V(\G)$.
\end{lem}

\f\demo Since $t \in \mz_{m}^*$ by Definition \ref{multi}, for each $i\in\mz_n$, $\G[V_i] \cong \C_{s}\circ \frac{m}{s}\K_1$. It suffices to prove that $V_0$ is a block of imprimitivity of $\Aut(\G)$ on $V(\G)$. Suppose that $V_0^g\cap V_0\neq \emptyset$ for some $g\in\Aut(\G)$. Take $v_0=u_0^g\in V_0^g\cap V_0$. Suppose that $v_0$ has a neighbour, say $x_0^g$ in $V_0^g$ but not in $V_0$. Since the edges between $V_i$ and $V_{i+1}$ are independent for any $i\in\mz_n$, $v_0$ is the only neighbour of $x_0^g$ in $V_0$.

Since $\G[V_0] \cong \C_{s}\circ \frac{m}{s}\K_1$, we may assume that all the vertices in $\{u_0, u_1, \ldots, u_{\frac{m}{s}-1}\}$ have the same neighbourhood in $\G[V_0]$. This implies that all vertices in the set $U=\{v_0=u_0^g, u_1^g, \ldots, u_{\frac{m}{s}-1}^g\}$ have the same neighbourhood in $\G[V_0^g]$. In particular, $x_0^g$ is adjacent to each vertex in $U$. Since $v_0$ is the only neighbour of $x_0^g$ in $V_0$, we have $U \cap V_0=\{v_0\}$. Note that $v_0$ has only two neighbours outside of $V_0$, and $\G$ has valency $2\left(\frac{m}{s}+1\right)$. Hence $|\G(v_0)\cap V_0\cap V_0^g|\geq 2\left(\frac{m}{s}-1\right) \geq 2$ as $\frac{m}{s}\geq 2$. Observe that each vertex in $U$ is adjacent to all vertices in $ \G(v_0)\cap V_0\cap V_0^g$. However, since $u_1^g\notin V_0$, $u_1^g$ has at most one neighbour in $V_0$, which is a contradiction. Thus all neighbours of $v_0$ in $V_0^g$ are contained in $V_0$. By the arbitrariness of $v_0$, we have $V_0^g\subseteq V_0$ and so $V_0^g=V_0$, completing the proof.\hfill\qed

\section{Proof of Theorem \ref{exist}}
\label{sec:exist}

The purpose of this section is to prove the following result, which together with Lemma~\ref{cay-meta} implies Theorem~\ref{exist}.

\begin{theorem}
\label{th}
Let $p$ be an odd prime, $m$ and $n$ be integers with $m\geq n+2 \geq 3$, and $\ld$ be an element of $\mz_{p^m}^*$ with order $p^{n+1}$. Then ${\rm MP}_{p^m, p^n, p^{m-1}, \ld}$ is a weak metacirculant which is a Cayley graph but not a weak metacirculant Cayley graph.
\end{theorem}

In the rest of this section, we always let $p, m, n$ and $\ld$ be as in Theorem \ref{th}, and $H=\lg h\rg\cong C_m$ be as in Definition \ref{multi}.
By Definition~\ref{multi}, ${\rm MP}_{p^m, p^n, p^{m-1}, \ld}$ has vertex set $\cup_{i=0}^{p^n-1}V_i$ and edge set $\cup_{i=0}^{p^n-1}(E_i \cup E_{i, i+1})$, reading the subscripts modulo $p^n$, where
$$
\begin{array}{lll}
V_i&=&\{(h^j, i)\ |\ j\in\mz_{p^m}\},\\
E_i&=&\{\{(h^j, i), (h^{j+kp^{m-1}+\ld^i}, i)\}, \{(h^j, i), (h^{j-kp^{m-1}-\ld^i}, i)\}\ |\ k\in\mz_p, j\in\mz_{p^m}\},\\
E_{i, i+1}&=&\{\{(h^j, i), (h^j, {i+1})\}\ |\ j\in\mz_{p^m}\}.
\end{array}
$$

The proof of Theorem~\ref{th} consists of the following three lemmas.

\begin{lem}\label{ismeta}
The graph ${\rm MP}_{p^m, p^n, p^{m-1}, \ld}$ is a metacirculant. Moreover, $\Aut({\rm MP}_{p^m, p^n, p^{m-1}, \ld})$ is transitive on the set of those arcs of ${\rm MP}_{p^m, p^n, p^{m-1}, \ld}$ whose underlying edges are in $\cup_{i=0}^{p^n-1} E_i$.
\end{lem}

\f\demo Denote $\G={\rm MP}_{p^m, p^n, p^{m-1}, \ld}$. Recall that for $g\in H$, $R(g)$ is the permutation on the vertices of $\G$ defined by:
$$
(h^j, i)^{R(g)}=(h^j g, i),\; \mbox{for}\; i\in\mz_{p^n}, j\in \mz_{p^m}.
$$
Recall also that $R(H)=\{R(g)\ |\ g\in H\}$ is a semiregular subgroup of $\Aut(\G)$ isomorphic to $H$ whose orbits on $V(\G)$ are $V_i$, $i \in \mz_{p^n}$.

Let $\a$ be the automorphism of $H$ such that $h^\a=h^{\ld}$. Define a permutation $\s_\a$ on the vertices of $\G$ by
$$(h^j, i)^{\s_\a}=((h^j)^\a, i+1),\; \mbox{for}\; i\in\mz_{p^n}, j\in \mz_{p^n}.$$
For each $i\in\mz_{p^n}$, it is easy to see that $E_{i, i+1}^{\s_\a}=E_{i+1, i+2}$.
Furthermore,
$$\begin{array}{lll}
\{(h^j, i), (h^{j+kp^{m-1}+\ld^i}, i)\}^{\s_\a}&=&\{(h^{j\ld}, {i+1}), (h^{j\ld+k\ld p^{m-1}+\ld^{i+1}}, {i+1})\}\in E_{i+1},\\
\{(h^j, i), (h^{j-kp^{m-1}-\ld^i}, i)\}^{\s_\a}&=&\{(h^{j\ld}, {i+1}), (h^{j\ld-k\ld p^{m-1}-\ld^{i+1}}, {i+1})\}\in E_{i+1},
\end{array}
$$
and so $E_{i}^{\s_\a}=E_{i+1}$ for each $i\in\mz_{p^n}$. This implies that $\s_\a$ preserves the adjacency relation of $\G$, and so $\s_\a\in\Aut(\G)$.

For any $(h^j, i)\in V(\G)$, we have $(h^j, i)^{R(h)\s_\a}=(h^{(j+1)\ld}, i+1)=(h^j, i)^{\s_\a R(h^\ld)}$. It follows that $R(h)^{\s_\a}=R(h^{\ld})=R(h)^{\ld}$, and so $\lg R(h), \s_\a\rg$ is metacyclic.
Clearly, $\lg R(h), \s_\a\rg$ is transitive on $V(\G)$, and $((1, 0), (1, 1), \ldots, (1, {p^n-1}))$ is a cycle of $\s_\a$ (as a permutation on $V(\G)$). So $\G$ is a metacirculant.

To show the second statement, we first observe that $\Aut(\G)$ preserves the edges in $\cup_{i=0}^{p^n-1} E_i$. The subgraph of $\G$ induced by $V_0$ is $\G_0=(V_0, E_0)$. It can be verified that $\lg R(h)\rg: \lg\s_\a^{p^{n}}\rg\cong C_{p^m}: C_p$ acts transitively on $E_0$. Let $\b$ be the automorphism of $H$ inverting every element of $H$. Let $\s_\b$ be a permutation of $V(\G)$ such that
$$
(h^j, i)^{\s_\b}=((h^j)^\b, i),\; \mbox{for}\; i\in\mz_{p^n}, j\in \mz_{p^n}.
$$
One can verify that $\s_\b\in\Aut(\G)$ and $\s_\b$ fixes each $V_i$ setwise. Furthermore, ${R(h^{-1})\s_\b}$ takes the arc $((1, 0), (h, 0))$ to its inverse arc $((h, 0), (1, 0))$. This implies that $\lg R(h), \s_\a^{p^{n}}, \s_\b\rg$ is transitive on the set of arcs of $\G_0$. Since $\s_\a$ cyclically permutes $V_i$'s, it follows that $\Aut(\G)$ is transitive on those arcs of $\G$ whose underlying edges are in $\cup_{i=0}^{p^n-1} E_i$.
\hfill\qed

Let
$$
\mathcal{G}=\lg x, y, z\ |\ x^{p^{m-1}}=y^{p^n}=z^p=1,\ y^{-1}xy=x^\ld,\ [z, x]=[z, y]=1\rg.
$$
It is easily seen that $\lg x^{p^{m-2}}, y^{p^{n-1}}, z\rg\cong C_p^3$, and so $\mathcal{G}$ is a non-metacyclic group.

\begin{lem}
${\rm MP}_{p^m, p^n, p^{m-1}, \ld} \cong \Cay(\MG, S\cup S^{-1})$, where
$$
S=\{x, xz, xz^2, \ldots, xz^{(p-1)}, y\}.
$$
\end{lem}

\demo Denote $\G={\rm MP}_{p^m, p^n, p^{m-1}, \ld}$ and $\Sigma'=\Cay(\MG, S\cup S^{-1})$. Define
$$
f: y^ix^{j}z^k\mapsto (h^{kp^{m-1}+j}, {i}),\; \mbox{for}\; i\in\mz_{p^n}, j\in\mz_{p^{m-1}}, k\in\mz_{p}.
$$
It can be verified that $f$ is a bijection from $V(\Sigma')$ to $V(\G)$. The neighbourhood of $y^ix^{j}z^k$ in $\Sigma'$ is
$$
\{y^i x^{j \pm \ld^i}z^{k \pm l}\ |\ l\in\mz_p\} \cup \{y^{i+1}x^{j}z^k, y^{i-1}x^jz^k\}.
$$
The image of this set under $f$ is
$$
\{(h^{kp^{m-1}+j \pm (lp^{m-1}+\ld^i)},i)\ |\ l\in\mz_p\}\cup\{(h^{kp^{m-1}+j}, {i+1}), (h^{kp^{m-1}+j}, {i-1})\},
$$
which is exactly the neighbourhood of $f(y^ix^{j}z^k)=(h^{kp^{m-1}+j}, i)$ in $\G$. Therefore, $f$ is an isomorphism from $\Sigma'$ to $\G$.
\hfill\qed

\begin{lem}\label{non-meta-cay}
The graph ${\rm MP}_{p^m, p^n, p^{m-1}, \ld}$ is not a weak metacirculant Cayley graph.
\end{lem}

\f\demo Denote $\G={\rm MP}_{p^m, p^n, p^{m-1}, \ld}$ and $A=\Aut(\G)$. We first prove the following claim.\medskip

\medskip
\noindent {\bf Claim.}\; For any $i\in\mz_{p^n}$ and $j\in\mz_{p^m}$, if $j\not\equiv 0\ (\mod p^{m-1})$, then the distance between $(1, i)$ and $(h^j, i)$ in $\G[V_i]$ is $\min\{t, p^{m-1}-t\}$, where $t\ld^i\equiv j\ (\mod p^{m-1})$ and $0\leq t\leq  p^{m-1}-1$.\medskip

Given $i\in\mz_{p^n}$ and $0\leq\ell\leq p^{m-1}-1$, let
$$
V_{i\ell}=\{(h^{kp^{m-1}+\ell}, i)\ |\ k\in\mz_p\}.
$$
Then
$$\B=\{V_{i\ell}\ |\ 0\leq\ell\leq p^{m-1}-1\}$$
is a partition of $V_i$. Moreover, each $V_{i\ell}$ is an independent set of $\G$, and the subgraph of $\G$ induced by $V_{i\ell}\cup V_{i(\ell+\ld^i)}$ is isomorphic to $\K_{p,p}$. Since $\G[V_i] \cong \C_{p^{m-1}}\circ p\K_1$, the quotient graph $Y$ of $\G[V_i]$ relative to $\B$ is a cycle of length $p^{m-1}$.

For any $j\in\mz_{p^m}$, if $j\not\equiv 0\ (\mod p^{m-1})$, then there exists $1\leq\ell\leq p^{m-1}-1$ such that $j\equiv \ell\ (\mod p^{m-1})$, and moreover, there exists $0\leq t\leq p^{m-1}-1$ such that $t\ld^i\equiv \ell\ (\mod p^{m-1})$. The distance $d$ between $(1, i)$ and $(h^j, i)$ in $\G[V_i]$ is just the distance between $V_{i0}$ and $V_{i\ell}=V_{i(t\ld^i)}$ in the quotient graph $Y$. It follows that $d=\min\{t, p^{m-1}-t\}$, completing the proof of the Claim.\medskip

By Lemma~\ref{block}, for each $i\in\mz_{p^n}$, $V_i$ is a block of imprimitivity of $A$ on $V(\G)$. It follows that $\cup_{i=0}^{p^n-1}\G[V_i]$ is a $2p$-factor of $\G$ which is invariant under the action of $A$.  So $A$ preserves the set
$$
F=\cup_{i=0}^{p^n-1} E_i.
$$
Let $\G'=\G-F$. Then $A\leq\Aut(\G')$, and for each $j\in\mz_{p^m}$, $\G[B_j]$ is a subgraph of $\G'$ isomorphic to $\C_{p^n}$, where
$$
B_j=\{(h^j, 0), (h^j, 1), (h^j, 2), \ldots, (h^j, {p^n-1})\}.
$$
So $\G' \cong p^m\C_{p^n}$ is the union of $p^m$ vertex-disjoint cycles of length $p^n$.

Suppose that $\G$ is a Cayley graph of a metacyclic $p$-group $G$. Say, $\G \cong \Lambda=\Cay(G, S)$ for an inverse-closed subset $S$ of $G \setminus \{1_G\}$ that generates $G$. Recall that $\G'\cong p^m\C_{p^n}$ is a $2$-factor of $\G$ invariant under $A$, $\G-E(\G')\cong p^n(\C_{p^{m-1}}\circ p\K_1)$, and for any two components of $\G-E(\G')$, either there is no edge connecting them in $\G$ or the edges between them form a perfect matching. Since $\G \cong \Lambda$, $\Lambda$ also has a $2$-factor $\Lambda' \cong p^m\C_{p^n}$ invariant under $\Aut(\Lambda)$ such that $\Lambda-E(\Lambda')\cong p^n(\C_{p^{m-1}}\circ p\K_1)$ and for any two components of $\Lambda-E(\Lambda')$, either there is no edge connecting them in $\Lambda$ or the edges between them form a perfect matching. So there exists $x\in S$ such that $\{1_G, x\}$ is an edge of $\Lambda'$. Since $\Aut(\Lambda)\leq\Aut(\Lambda')$, $\{1_G, x\}^{R(x^\ell)}=\{x^\ell, x^{\ell+1}\}\in E(\Lambda')$ for any $0\leq\ell\leq o(x)$, which implies that $(1_G, x, x^2, \ldots, x^{o(x)-1})$ is cycle of $\Lambda'$. Since $\Lambda'\cong p^m\C_{p^n}$, we have $o(x)=p^n$. Clearly, $x^{-1}\in S$ and $\Cay(G, \{x, x^{-1}\})\cong p^m\C_{p^n}$. Recall that $G$ acts faithfully on $V(\Lambda)$ by right multiplication, and this action induces a regular subgroup $R(G)$ of $\Aut(\Lambda)$, which will be identified with $G$ in the sequel. Let $P$ be a Sylow $p$-subgroup of $\Aut(\Lambda)$ such that $G\leq P$. From the proof of Lemma~\ref{ismeta}, we see that $G$ is a proper subgroup of $P$, and so $N_P(G)>G$. It follows that $p\ |\ |\Aut(G, S)|$. Let $\b\in \Aut(G, S)$ be of order $p$. Then $\b$ preserves $\Cay(G, \{x, x^{-1}\})\cong p^m\C_{p^n}$, and so $\b$ must fix $x$. Consequently, $p\ |\ |\Aut(G, S-\{x, x^{-1}\})|$, and hence $S=\{x, x^{-1}\}\cup \{y^{\b^i}, (y^{-1})^{\b^i}\ |\ i\in\mz_p\}$ for some $y\in S-\{x,x^{-1}\}$.

Since $\Lambda-E(\Lambda')\cong p^n (\C_{p^{m-1}}\circ p\K_1)$, we have  $\Lambda-E(\Lambda')=\Cay(G, S-\{x, x^{-1}\})$. Let $T=\lg S-\{x, x^{-1}\}\rg$. Then $\Cay(T, S-\{x,x^{-1}\})$ is a subgraph of $\Lambda$ isomorphic to $\C_{p^{m-1}}\circ p\K_1$. So $T$ has order $p^m$. Then $G=\cup_{i=0}^{p^n-1}Tx^i$. Recall that the $p^m$ edges of $\Lambda$ between $T$ and $Tx$ are independent. So the quotient graph of $\Lambda$ relative to $\{Tx^i\ |\ i\in\mz_{p^n}\}$ is a cycle of length $p^n$. Let $K$ be the kernel of $G$ acting on $\{Tx_i\ |\ i\in\mz_{p^n}\}$. Then $G/K\leq D_{2p^n}$ and so $G/K\cong\mz_{p^n}$ as $G$ is a $p$-group. This implies that $T=K$. We claim that $T$ is cyclic. Suppose on the contrary that $T$ is non-cyclic. Then $G$ is non-cyclic. Since $G$ is metacyclic and $p>2$, by Lemma~\ref{p2-group} the subgroup $\Omega_1(G)$ of $G$ generated by the elements of order $p$ is an elementary abelian group of order $p^2$. Hence $\lg x\rg\cap\Omega_1(G)\neq 1_G$ and $\Omega_1(G)\leq T$. Consequently, $T\cap \lg x\rg\neq 1_G$. Note that $\lg x\rg$ acts transitively on the $p^n$ components of $\Cay(G, S-\{x, x^{-1}\})$, and $T$ is the stabilizer of the block $T$ in $G$. Hence $G=T\lg x\rg$. Since $|T|=p^m$ and $o(x)=p^n$, from $|G|=p^{m+n}$ it follows that $\lg x\rg\cap T=1_G$, a contradiction. Thus, $T$ is a normal cyclic subgroup of $G$ of order $p^m$, and moreover, $G=T: \lg x\rg$. Therefore, $T=\lg y\rg$ and $S-\{x, x^{-1}\}=\{y^{kp^{m-1}+1}, y^{-kp^{m-1}-1}\ |\ k\in\mz_p\}$.

Assume that $y^x=y^{\ld'}$ for some $\ld'\in\mz_{p^m}^*$. Since $x$ has order $p^n$, $\ld'$ has order at most $p^n$. 
Recall that the edges of $\Lambda$ between $T$ and $Tx^{-1}$ form a perfect matching. Note that $1_G\sim x^{-1}$ and $y\sim x^{-1}y$. Since $x^{-1}yx=y^{\ld'}$, we have $x^{-1}y=y^{\ld'}x^{-1}$. We now consider the distance $d'$ between $x^{-1}$ and $x^{-1}y=y^{\ld'}x^{-1}$ in the subgraph induced by $Tx^{-1}$. Indeed, $d'$ is just the distance between $1_G$ and $y^{\ld'}$ in the subgraph induced by $T$. Observe that the subgraph induced by $T$ is the Cayley graph
$$
\Cay(T, \{y^{kp^{m-1}+1}, y^{-kp^{m-1}-1}\ |\ k\in\mz_p\}),
$$
which is isomorphic to $\C_{p^{m-1}}\circ p\K_1$. Let $B=\lg y^{p^{m-1}}\rg$. Then
$$
T=\cup_{\ell=0}^{p^{m-1}-1}By^\ell.
$$
It is clear that each coset $By^\ell$ is an independent set of $\Lambda[T]$, and the subgraph of $\Lambda$ induced by $By^{\ell}\cup By^{\ell+1}$ is isomorphic to $\K_{p, p}$. Note that $By^{\ld'}=By^t$, where $\ld'\equiv t'\ (\mod p^{m-1})$ and $1\leq t'\leq p^{m-1}-1$. Hence $d'=\min\{t', p^{m-1}-t'\}$.

Suppose that $f$ is an isomorphism from $\Lambda$ to $\Gamma$. Since $\Gamma$ is vertex-transitive, we may assume that $f$ maps $1_G$ to $(1, 0)$. By Lemma~\ref{ismeta}, the arcs in $\G[V_0]$ are equivalent under $A$. So we may further assume that $f$ takes the arc $(1_G, y)$ of $\Lambda$ to the arc $((1, 0), (h, 0))$ of $\G$. Clearly, $f$ maps $Tx^{-1}$ to $V_1$ or $V_{p^n-1}$.

If $f$ maps $Tx^{-1}$ to $V_1$, then since the edges between $T$ and $Tx$ form a perfect matching, $f$ maps $x^{-1}$ and $x^{-1}y$ ($=y^{\ld'}x^{-1}$) to $(1, 1)$ and $(h, 1)$, respectively. By the Claim above, the distance between $(1, 1)$ and $(h, 1)$ in $\G[V_1]$ is $\min\{t, p^{m-1}-t\}$, where $t\ld\equiv 1\ (\mod p^{m-1})$ and $1\leq t\leq p^{m-1}-1$. Hence $\min\{t', p^{m-1}-t'\} = \min\{t, p^{m-1}-t\}$.
It follows that $\ld^{-1}\equiv \ld'\ (\mod p^{m-1})$ because $\ld$ and $\ld'$ have odd orders, where $\ld^{-1}$ is the inverse of $\ld$ in $\mz_n^*$. Consequently, $\ld^{-1}=kp^{m-1}+\ld'$ for some $k\in\mz_{p^m}$. Hence $\ld^{-p}\equiv (\ld')^p\ (\mod p^m)$. This implies that the order of $\ld^{-1}\in\mz_{p^m}^*$ is at most $p^n$, a contradiction.

Similarly, if $f$ maps $Tx^{-1}$ to $V_{p^{n}-1}$, then $f$ maps $x^{-1}$ and $x^{-1}y$ ($=y^{\ld'}x^{-1}$) to $(1, {p^{n}}-1)$ and $(h, {p^{n}}-1)$, respectively. Again by the Claim, the distance between $(1, 1)$ and $(h, {p^{n}}-1)$ is $\min\{t, p^{m-1}-t\}$, where $t\ld^{p^{n}-1}\equiv 1\ (\mod p^{m-1})$ and $1\leq t\leq p^{m-1}-1$. Hence $\min\{t', p^{m-1}-t'\}=\min\{t, p^{m-1}-t\}$. It follows that $\ld^{1-p^{n}}\equiv \ld'\ (\mod p^{m-1})$ because $\ld$ and $\ld'$ have odd orders. Consequently, $\ld^{1-p^{n}}=kp^{m-1}+\ld'$ for some $k\in\mz_{p^m}$. Hence $(\ld^{1-p^{n}})^p\equiv (\ld')^p\ (\mod p^m)$. This implies that the order of $\ld^{1-p^{n}}\in\mz_{p^m}^*$ is at most $p^n$. However, $\ld^{1-p^{n}}$ and $\ld$ have the same order which is assumed to be $p^{n+1}$, a contradiction.
\hfill\qed

So far we have proved Theorem \ref{th}. As mentioned earlier, Theorem \ref{exist} follows from Theorem \ref{th} and Lemma \ref{cay-meta}.




\section{Proof of Theorem~\ref{2p+2}}
\label{sec:2p+2}

Let $p$ be an odd prime and $\G$ a connected metacirculant graph of order $p^4$ and valency $2p+2$.
By Theorem~\ref{p-metacirculant}, $\G$ has a split metacyclic vertex-transitive group of automorphisms, say $G$. We may further assume that $G$ is a $p$-group. If $Z(G)$ is non-cyclic, then $G$ is regular on $V(\G)$ by Lemma~\ref{cyclic-center}, and so $\G$ is a metacirculant Cayley graph. If $G$ is abelian, then again $G$ is regular on $V(\G)$ and $\G$ is a metacirculant Cayley graph.

In what follows, we assume that $Z(G)$ is cyclic and $G$ is non-abelian. From the proof of Theorem~\ref{p-metacirculant}, we may assume that $G=\lg x\rg: \lg y\rg\cong C_{p^m}: C_{p^n}$ for some $m\geq n$ and $G_v\leq \lg y\rg$ for some $v\in V(\G)$. Since $G_v$ is core-free in $G$, every non-identity element of $\lg y\rg$ induces a non-trivial automorphism of $\lg x\rg$ by conjugation. Since $\Aut(C_{p^m})\cong C_{p^{m-1}}: C_{p-1}$, it follows that $n\leq m-1$. Since $\lg x\rg\unlhd G$, the transitivity of $G$ on $V(\G)$ implies that $\lg x\rg$ acts semiregularly on $V(\G)$. Since $|V(\G)|=p^4$, we have $m\leq 4$ and $m+n \geq 4$. If $m=4$, then $\lg x\rg$ acts regularly on $V(\G)$ and so $\G$ is a metacirculant Cayley graph. If $m<4$ and $m+n=4$, then $G$ acts regularly on $V(\G)$ and hence $\G$ is a metacirculant Cayley graph.

Assume $m < 4 < m+n$. Then the only possibility is $(m, n)=(3, 2)$, which implies that $G_v=\lg y^p\rg\cong C_p$ and $y^{-1}xy=x^{kp+1}$ for some $k\in\mz_p^*$. Consequently, $G' = \lg x^{p}\rg$ ($\cong C_{p^2}$).
Since $G=\lg x\rg: \lg y\rg\cong C_{p^3}: C_{p^2}$ and $p>2$, from Lemma~\ref{p2-group} it follows that $\Omega_1(G)=\lg x^{p^{2}}\rg\times\lg y^{p}\rg\cong C_p\times C_p$. By the N/C theorem, we have $G/C_G(\Omega_1(G))\leq\Aut(\Omega_1(G))\cong {\rm GL}(2, p)$. Since $|{\rm GL(2, p)}|=(p^2-p)(p^2-1)$ and $G$ is a $p$-group, we have $G/C_G(\Omega_1(G))=1$ or $G/C_G(\Omega_1(G))\cong C_p$. The former cannot happen, for otherwise $\Omega_1(G)$ is contained in $Z(G)$ and hence $G_v\unlhd G$, a contradiction. So $G/C_G(\Omega_1(G))\cong C_p$.  Since $G_v\leq\Omega_1(G)$, $C_G(\Omega_1(G))\leq C_G(G_v)\leq N_G(G_v)$ and hence $C_G(\Omega_1(G))=C_G(G_v)=N_G(G_v)$ as $G_v$ is non-normal in $G$.

By Proposition~\ref{coset}, $\G$ is isomorphic to the coset graph $\G'=\Cos(G, G_v, G_vSG_v)$, where $S$ consists of the elements of $G$ each of which maps $v$ to one of its neighbours. We will simply identify $\G$ with $\G'$ in the remainder of the proof. Since $\G$ is connected, we have $G=\lg G_vSG_v\rg$, which implies that there exists $d\in S\setminus C_G(G_v)$. Then $d=x^iy^j$ for some $i\in\mz_{p^3}^*$ and $j\in\mz_{p^2}$, and furthermore, $G_v^d\neq G_v$ and so $G_v^d\cap G_v=1_G$ as $G_v\cong C_p$. Consequently, $|G_vdG_v|/|G_v|=p$. If $d^{-1}\in G_vdG_v$, then the subgraph $\Cos(\lg G_vdG_v\rg, G_v, G_vdG_v)$ of $\G$ would have odd order and odd  valency $p$, but this cannot happen. Thus $d^{-1}\notin G_vdG_v$ and $|G_v\{d, d^{-1}\}G_v|/|G_v|=2p$ by Lemma~\ref{lem-coset}. We may assume that $D=G_v\{d, d^{-1}\}G_v\cup G_v\{c, c^{-1}\}G_v$ for some $c\in G$. Since $\G$ has valency $2p+2$, we have $|G_v\{c, c^{-1}\}G_v|/|G_v|=2$. It follows that $c$ normalizes $G_v$ and so $c\in N_G(G_v)=C_G(G_v)$.

Let $M=\lg d, G_v\rg$. As $d\notin C_G(G_v)=\lg x^p\rg: \lg y\rg$, $o(d)=p^3$ and so $\lg d^{p^2}\rg=\lg x^{p^2}\rg$. It follows that  $\Omega_1(G)\leq M$ and $M/\Omega_1(G)=\lg d\Omega_1(G)\rg\cong C_{p^2}$. Consequently, $|M|=p^4$ and so $M\unlhd G$.

Let $\Sigma=\Cos(M, G_v, G_v\{d, d^{-1}\}G_v)$. Then $\Sigma$ has order $p^3$ and valency $2p$, and $M$ is a vertex- and edge-transitive group of automorphisms of $\G_1$ (see Lemma~\ref{lem-coset}). Recall that $\Omega_{1}(G)=\lg d^{p^2}\rg\times G_v\cong C_p\times C_p$. The quotient graph $\Sigma_{\Omega_1(G)}$ of $\Sigma$ relative to $\Omega_1(G)$ has $p^2$ vertices, $\Omega_1(G)$ is the kernel of $M$ acting on $V(\Sigma_{\Omega_1(G)})$, and $M/\Omega_1(G)$ is edge-transitive on $\Sigma_{\Omega_1(G)}$. Since $G_v\leq\Omega_1(G)$, the subgraph induced by some pair of adjacent orbits of $\Omega_1(G)$ is isomorphic to $\K_{p, p}$. Since $M/\Omega_1(G)$ is edge-transitive on $\Sigma_{\Omega_1(G)}$, it follows that the subgraph induced by any two adjacent orbits of $\Omega_1(G)$ is isomorphic to $\K_{p, p}$. Consequently, $\Sigma_{\Omega_1(G)}\cong \C_{p^2}$ and $\Sigma \cong \C_{p^2}\circ p\K_1$. Since $G$ acts transitively on $V(\G)$, this implies that the subgraph of $\G$ induced by each orbit of $M$ is isomorphic to $\C_{p^2}\circ p\K_1$.

Consider the quotient graph $\G_M$ of $\G$ relative to $M$. Since $M\unlhd G$ and the subgraph of $\G$ induced by each orbit of $M$ is isomorphic to $\C_{p^2}\circ p\K_1$, we have $\G_M\cong {\C}_p$. So we may assume $V(\G_M)=\{B_i\ |\ i\in \mz_p\}$ such that $B_i\sim B_{i+1}$ for $i\in\mz_p$. Then $\G[B_i\cup B_{i+1}] \cong p^3\K_2$ for each $i\in\mz_p$. By a similar argument as in the proof of Lemma~\ref{block}, one can show that $V(\G_M)$ is a system of blocks of imprimitivity of $\Aut(\G)$.

Assume that $c^p\notin G_v$. Since $c\in C_{G}(\Omega_1(G))=N_G(G_v)\cong C_{p^2}: C_{p^2}$, $c$ has order $p^2$ and so $(G_v, G_vc, \ldots, G_vc^{p^2-1})$ is a cycle of $\G$ of length $p^2$. Furthermore, $\lg c^p\rg\times G_v=\Omega_1(G)$, so $\{G_v, G_vc^p, G_vc^{2p}, \ldots, G_vc^{(p-1)p}\}$ is an orbit of $\Omega_1(G)$, and the vertices in this orbit have the same neighbourhood in $\G_1$.

If $\G$ is not a Cayley graph, then case (b) in Theorem \ref{2p+2} holds. In the sequel we assume that $\G$ is a Cayley graph of a group $N$, say $\G=\Cay(N, S)$. Let $K$ be the kernel of $N$ acting on $V(\G_M)$. Then $N/K\cong C_p$. Recall that $\G[B_0] \cong {\C}_{p^2}\circ p\K_1$. We call a subset of $B_0$ a {\em part} of $B_0$ if it is maximal with respect to the property that all vertices in the set have the same neighbourhood in $\G[B_0]$. Since $\G[B_0] \cong {\C}_{p^2}\circ p\K_1$, $B_0$ has $p^2$ parts, and each of them is a block of imprimitivity of $\Aut(\G[B_0])$ on $B_0$. Let $L$ be the kernel of $K$ acting on the $p^2$ parts of $\G[B_0]$. Then $K/L\cong C_{p^2}$ and $L\cong C_p$. This implies that $K\cong C_{p^3}$ or $K\cong C_{p^2}\times C_{p}$. In the former case, $N$ is a metacyclic group, and so $\G$ is a metacirculant Cayley graph. Suppose $K\cong C_{p^2}\times C_{p}$ in what follows.

Since $V(\G_M)$ is a system of blocks of imprimitivity of $\Aut(\G)$ on $V(\G)$, we have  $\Aut(\G)\leq\Aut(\G^*)$, where $\G^*$ is obtained from $\G$ by deleting the edges contained in each $\G[B_i]$. Since $(G_v, G_vc, \ldots, G_vc^{p^2-1})$ is a cycle of length $p$, we have $\G^*\cong p^2\C_{p^2}$. Since $\G=\Cay(N, S)$, we may relabel the vertices of $\G$ by the elements of $N$. Then $S$ contains an element $g\in N\setminus K$ such that $o(g)=p^2$ and $(1_N, g, g^2, \ldots, g^{p^2-1})$ is a cycle of $\G$ corresponding to $(G_v, G_vc, \ldots, G_vc^{p^2-1})$. Recall that $\{G_v, G_vc^p, G_vc^{2p}, \ldots, G_vc^{(p-1)p}\}$ is an orbit of $\Omega_1(G)$ that is also a part of $B_0$. We can label this part by $\{1, g^p, g^{2p}, \ldots, g^{(p-1)p}\}$. Note that $N$ acts on $V(G)=N$ by right multiplication. So $\lg g^p\rg$ fixes each of the $p^2$ parts of $B_0$. It then follows that $K/\lg g^p\rg\cong C_{p^2}$. So we may assume that $K=\lg g^p\rg\times\lg g_1\rg\cong C_p\times C_{p^2}$. Then $N=\lg g\rg\lg g_1\rg$, and $N$ is metacyclic as $p>2$. (Note that, by \cite[III, 11.5]{Huppert}, if $G=\lg a\rg\lg b\rg$ is a $p$-group with $p>2$, then $G$ is metacyclic.) This implies that $\G$ is a metacirculant Cayley graph.

Now assume that $c^p\in G_v$. Then $c$ has order $p^2$. Assume that $v\in B_0$. We may label the vertices in $B_0$ in the following way: $v=(d^0, 0)$ and $(d^i, 0)=v^{G_vd^i}$ for $i\in\mz_{p^3}$. Note that $\G[B_0] \cong {\C}_{p^2}\circ p\K_1$ and $\Omega_1(G)$ is the kernel of $M$ acting on the $p^2$ parts of $\C_{p^2}\circ p\K_1$. As $M=G_v\lg d\rg$, we have $(d^0, 0)^{\Omega_1(G)}=\{(d^{kp^2}, 0)\ |\ k\in\mz_p \}$. Since $v$ is adjacent to $v^d$, $(d^0, 0)$ is adjacent to $(d^1, 0)$. Hence $(d^0, 0)$ is adjacent to all vertices in $(d^1, 0)^{\Omega_{1}(G)}=\{(d^{1+kp^2}, 0)\ |\ k\in\mz_p\}$. Similarly, since $v$ is adjacent to $v^{d^{-1}}$, $(d^0, 0)$ is adjacent to all vertices in $(d^{-1}, 0)^{\Omega_{1}(G)}=\{(d^{-1-kp^2}, 0)\ |\ k\in\mz_p\}$. Now the edge set of $\G[B_0]$ is $E_0=\{\{(d^i, 0), (d^{i+1+kp^2}, 0)\}, \{(d^i, 0), (d^{i-1-kp^2}, 0)\}\ |\ i\in\mz_{p^3}, k\in\mz_p\}$.

Since $c$ cyclically permutates the $p$ orbits of $M$, we may assume that $B_j=B_0^{c^j}$ for $j\in\mz_p$. Since $v$ is adjacent $v^c$,
$(v, v^c, v^{c^2}, \ldots, v^{c^{p-1}})$ is a cycle of length $p$. Without loss of generality we may assume $(d^0, j)=v^{c^j}$ for $j\in\mz_p$ and $(d^i, j)=(d^0, j)^{d^i}$ for $i\in\mz_{p^3}$. We now have $B_j=\{(d^i, j)\ |\ i\in\mz_{p^3}\}$ for $j\in\mz_p$ and $(d^i, j)\sim (d^i, j+1)$ for $i\in\mz_{p^3}, j\in\mz_p$.

Since $(d^0, 0)^{c^j}=(d^0, j)$, the set of neighbours of $(d^0, j)$ in $B_j$ consists of the two orbits of $\Omega_1(G)$ containing $(d, 0)^{c^j}$ and $(d^{-1}, 0)^{c^j}$ respectively, namely,
\begin{equation}\label{eq-2}
((d^0, 0)^{dc^j})^{\Omega_{1}(G)}\cup ((d^0, 0)^{d^{-1}c^j})^{\Omega_{1}(G)}=((d^0, 0)^{c^j})^{c^{-j}dc^j\Omega_{1}(G)}\cup ((d^0, 0)^{c^j})^{c^{-j}d^{-1}c^j\Omega_{1}(G)}.
\end{equation}
Since $M/\Omega_1(G)$ is a normal cyclic subgroup of $G/\Omega_1(G)$ of order $p^2$, we have $c^{-1}dc\Omega_1(G)=d^{kp+1}\Omega_{1}(G)$ for some $k\in\mz_{p}$. Since $G'\nleqslant \Omega_{1}(G)$ (as $G'\cong C_{p^2}$), we have $k\in\mz_p^*$. Set $\ld=kp+1$. Then  $c^{-j}dc^j\Omega_1(G)=d^{\ld^j}\Omega_{1}(G)$. It follows that the right-hand side of equation \eqref{eq-2} is
$\{(d^{\ld^j+1+kp^2}, j), (d^{-\ld^j-1-kp^2}, j)\ |\ i\in\mz_{p^3}, k\in\mz_p\}$. So the edge set of $\G[B_j]$ is $$E_j=\{\{(d^i, j), (d^{i+\ld^i+kp^2}, j)\}, \{(d^i, j), (d^{i-\ld^i-kp^2}, j)\}\ |\ i\in\mz_{p^3}, k\in\mz_p\}.$$
Now we can see that $\G\cong{\rm MP}_{p^3, p^2, p^2,\ld}$. Since $\ld=kp+1$, $\ld$ must be an element of $\mz_{p^3}^*$ of order $p^2$. This completes the proof of Theorem~\ref{2p+2}.
\hfill\qed

\medskip
\f {\bf Acknowledgements}\; We appreciate the anonymous referees for their helpful comments. 
The first author was partially supported by the National
Natural Science Foundation of China (11671030, 11271012), the Fundamental
Research Funds for the Central Universities (2015JBM110) and the
111 project of China (B16002). The second author was supported
by the Australian Research Council (FT110100629).

{\small
	
{}
}


\begin{thebibliography}{99}

\bibitem{AP}
B. Alspach, T.D. Parson, A Construction for vertex-transitive
graphs, Canad. J. Math. 34 (1982) 307--318.

\bibitem{p-group1}
Y. Berkovich, Z. Janko, Groups of Prime Power Order, Volume 1, Walter de Gruyter GmbH \& Co. KG, Berlin, 2008.

\bibitem{p-group2}
Y. Berkovich, Z. Janko, Groups of Prime Power Order, Volume 2, Walter de Gruyter GmbH \& Co. KG, Berlin, 2008.

\bibitem{B}
N. Biggs, Algebraic Graph Theory, Second Edition, Cambridge
University Press, Cambridge, 1993.


\bibitem{CPZ}
M. Conder, T. Pisanski, A. \v Zitnik,
GI-graphs: a new class of graphs with many symmetries, J. Algebr. Comb. 40 (2014) 209--231.

\bibitem{Feng}
Y.-Q. Feng, On vertex-transitive graphs of odd prime-power
order, Discrete Math. 248 (2002) 265--269.

\bibitem{FengLuXu}
Y.-Q. Feng, Z.P. Lu, M.Y. Xu, Automorphism groups of Cayley
digraphs, in: Application of Group Theory to Combinatorics, edited by
J. Koolen, J.H. Kwak, M.Y. Xu, Taylor $\&$ Francis Group, London,
2008; pp. 13--25.


\bibitem{Godsil1981}
C.D. Godsil, On the full automorphism group of a graph, Combinatorica 1 (1981)
243--256.


\bibitem{Huppert}
B. Huppert, Eudiche Gruppen I, Springer-Verlag, 1967.

\bibitem{LLZ}
C.H. Li, Z.P. Lu, H. Zhang, Tetravalent edge-transitive Cayley graphs with odd number of vertices,
J. Combin. Theory Ser. B 96 (2006) 164--181.

\bibitem{LWS}
C.H. Li, S.J. Song, D.J. Wang, A characterization of metacirculants, J. Combin. Theory Ser. A 120 (2013) 39--48.


\bibitem{MS}
D. Maru\v si\v c, P. \v Sparl, On quartic half-arc-transitive metacirculants, J. Algebr. Comb.  28 (2008) 365--395.


\bibitem{Pan-thesis}
J.M. Pan, Groups with metacyclic factors, metacirculants and locally primitive graphs, Ph.D. Thesis, Yunnan University, 2009.

\bibitem{Redei}
Das L. R\' edei, ``schiefe Produkt" in der Gruppentheorie mit Anwendung auf die endlichen nichkommutativen Gruppen mit lauter
kommutative echten Untergruppen und Ordnungszahlen, zu denen nur kommutative Grouppen geh\" oren. Comment. Math. Helvet. 20 (1947) 225--264.


\bibitem{SPP}
M.L. Sara\v zin, Walter Pacco, A. Previtali, Generalizing the generalized Petersen graphs, Discrete Math. 307 (2007) 534--543.

\bibitem{Sabidussi}
B.O. Sabidussi, Vertex-transitive graphs, Monash. Math. 68 (1964)
426--438.

\bibitem{W}
M.E. Watkins,  A theorem on Tait colorings with an application to the generalized Petersen graphs, J. Combin. Theory 6 (1969) 152--164.

\bibitem{WI}
H. Wielandt, Finite Permutation Groups, Academic Press, New York,
1964.

\bibitem{X1}
M.Y. Xu, Automorphism groups and isomorphisms of Cayley digraphs,
Discrete Math. 182 (1998) 309--319.

\bibitem{Xu-Zhang}
M.Y. Xu, Q. Zhang, A classification of metacyclic $2$-groups, Algebra Colloq. 13 (2006) 25--34.

\end{thebibliography}
\end{document}